\magnification=1200

\hsize=11.25cm    
\vsize=18cm       
\parindent=12pt   \parskip=5pt     

\hoffset=.5cm   
\voffset=.8cm   

\pretolerance=500 \tolerance=1000  \brokenpenalty=5000

\catcode`\@=11

\font\eightrm=cmr8         \font\eighti=cmmi8
\font\eightsy=cmsy8        \font\eightbf=cmbx8
\font\eighttt=cmtt8        \font\eightit=cmti8
\font\eightsl=cmsl8        \font\sixrm=cmr6
\font\sixi=cmmi6           \font\sixsy=cmsy6
\font\sixbf=cmbx6

\font\tengoth=eufm10 
\font\eightgoth=eufm8  
\font\sevengoth=eufm7      
\font\sixgoth=eufm6        \font\fivegoth=eufm5

\skewchar\eighti='177 \skewchar\sixi='177
\skewchar\eightsy='60 \skewchar\sixsy='60

\newfam\gothfam           \newfam\bboardfam

\def\tenpoint{
  \textfont0=\tenrm \scriptfont0=\sevenrm \scriptscriptfont0=\fiverm
  \def\rm{\fam\z@\tenrm}
  \textfont1=\teni  \scriptfont1=\seveni  \scriptscriptfont1=\fivei
  \def\oldstyle{\fam\@ne\teni}\let\old=\oldstyle
  \textfont2=\tensy \scriptfont2=\sevensy \scriptscriptfont2=\fivesy
  \textfont\gothfam=\tengoth \scriptfont\gothfam=\sevengoth
  \scriptscriptfont\gothfam=\fivegoth
  \def\goth{\fam\gothfam\tengoth}
  
  \textfont\itfam=\tenit
  \def\it{\fam\itfam\tenit}
  \textfont\slfam=\tensl
  \def\sl{\fam\slfam\tensl}
  \textfont\bffam=\tenbf \scriptfont\bffam=\sevenbf
  \scriptscriptfont\bffam=\fivebf
  \def\bf{\fam\bffam\tenbf}
  \textfont\ttfam=\tentt
  \def\tt{\fam\ttfam\tentt}
  \abovedisplayskip=12pt plus 3pt minus 9pt
  \belowdisplayskip=\abovedisplayskip
  \abovedisplayshortskip=0pt plus 3pt
  \belowdisplayshortskip=4pt plus 3pt 
  \smallskipamount=3pt plus 1pt minus 1pt
  \medskipamount=6pt plus 2pt minus 2pt
  \bigskipamount=12pt plus 4pt minus 4pt
  \normalbaselineskip=12pt
  \setbox\strutbox=\hbox{\vrule height8.5pt depth3.5pt width0pt}
  \let\bigf@nt=\tenrm       \let\smallf@nt=\sevenrm
  \normalbaselines\rm}

\def\eightpoint{
  \textfont0=\eightrm \scriptfont0=\sixrm \scriptscriptfont0=\fiverm
  \def\rm{\fam\z@\eightrm}
  \textfont1=\eighti  \scriptfont1=\sixi  \scriptscriptfont1=\fivei
  \def\oldstyle{\fam\@ne\eighti}\let\old=\oldstyle
  \textfont2=\eightsy \scriptfont2=\sixsy \scriptscriptfont2=\fivesy
  \textfont\gothfam=\eightgoth \scriptfont\gothfam=\sixgoth
  \scriptscriptfont\gothfam=\fivegoth
  \def\goth{\fam\gothfam\eightgoth}
  
  \textfont\itfam=\eightit
  \def\it{\fam\itfam\eightit}
  \textfont\slfam=\eightsl
  \def\sl{\fam\slfam\eightsl}
  \textfont\bffam=\eightbf \scriptfont\bffam=\sixbf
  \scriptscriptfont\bffam=\fivebf
  \def\bf{\fam\bffam\eightbf}
  \textfont\ttfam=\eighttt
  \def\tt{\fam\ttfam\eighttt}
  \abovedisplayskip=9pt plus 3pt minus 9pt
  \belowdisplayskip=\abovedisplayskip
  \abovedisplayshortskip=0pt plus 3pt
  \belowdisplayshortskip=3pt plus 3pt 
  \smallskipamount=2pt plus 1pt minus 1pt
  \medskipamount=4pt plus 2pt minus 1pt
  \bigskipamount=9pt plus 3pt minus 3pt
  \normalbaselineskip=9pt
  \setbox\strutbox=\hbox{\vrule height7pt depth2pt width0pt}
  \let\bigf@nt=\eightrm     \let\smallf@nt=\sixrm
  \normalbaselines\rm}

\tenpoint

\def\pc#1{\bigf@nt#1\smallf@nt}         \def\pd#1 {{\pc#1} }

\catcode`\;=\active
\def;{\relax\ifhmode\ifdim\lastskip>\z@\unskip\fi
\kern\fontdimen2  -1.2 \fontdimen3 \string;}

\catcode`\:=\active
\def:{\relax\ifhmode\ifdim\lastskip>\z@\unskip\fi\penalty\@M\ \fi\string:}

\catcode`\!=\active
\def!{\relax\ifhmode\ifdim\lastskip>\z@
\unskip\fi\kern\fontdimen2  -1.1 \fontdimen3 \string!}

\catcode`\?=\active
\def?{\relax\ifhmode\ifdim\lastskip>\z@
\unskip\fi\kern\fontdimen2  -1.1 \fontdimen3 \string?}

\frenchspacing

\def\raggedbottom{\topskip 10pt plus 36pt\r@ggedbottomtrue}

\def\pointir{\unskip . --- \ignorespaces}

\def\Medbreak{\vskip-\lastskip\medbreak}

\long\def\th#1 #2\enonce#3\endth{
   \Medbreak\noindent
   {\pc#1} {#2\unskip}\pointir{\it #3}\smallskip}

\def\decale#1{\smallbreak\hskip 28pt\llap{#1}\kern 5pt}
\def\decaledecale#1{\smallbreak\hskip 34pt\llap{#1}\kern 5pt}
\def\puce{\smallbreak\hskip 6pt{$\scriptstyle\bullet$}\kern 5pt}

\def\eqalign#1{\null\,\vcenter{\openup\jot\m@th\ialign{
\strut\hfil$\displaystyle{##}$&$\displaystyle{{}##}$\hfil
&&\quad\strut\hfil$\displaystyle{##}$&$\displaystyle{{}##}$\hfil
\crcr#1\crcr}}\,}

\catcode`\@=12

\showboxbreadth=-1  \showboxdepth=-1

\newcount\numerodesection \numerodesection=1
\def\section#1{\bigbreak
 {\bf\number\numerodesection.\ \ #1}\nobreak\medskip
 \advance\numerodesection by1}

\mathcode`A="7041 \mathcode`B="7042 \mathcode`C="7043 \mathcode`D="7044
\mathcode`E="7045 \mathcode`F="7046 \mathcode`G="7047 \mathcode`H="7048
\mathcode`I="7049 \mathcode`J="704A \mathcode`K="704B \mathcode`L="704C
\mathcode`M="704D \mathcode`N="704E \mathcode`O="704F \mathcode`P="7050
\mathcode`Q="7051 \mathcode`R="7052 \mathcode`S="7053 \mathcode`T="7054
\mathcode`U="7055 \mathcode`V="7056 \mathcode`W="7057 \mathcode`X="7058
\mathcode`Y="7059 \mathcode`Z="705A


\def\hfl#1#2#3{\smash{\mathop{\hbox to#3{\rightarrowfill}}\limits
^{\textstyle#1}_{\textstyle#2}}}

\def\ogoth{{\goth o}}

\def\pgoth{{\goth p}}

\def\P{{\bf P}}
\def\Q{{\bf Q}}
\def\Qp{\Q_p}

\def\C{{\bf C}}
\def\N{{\bf N}}

\def\U{{\bf U}}
\def\Z{{\bf Z}}
\def\Zp{\Z_p}
\def\F{{\bf F}}
\def\Fp{{\F_{\!p}}}

\def\Aut{\mathop{\rm Aut}\nolimits}
\def\Frob{\mathop{\rm Frob}\nolimits}

\def\Jac{\mathop{\rm Jac}\nolimits}

\def\End{\mathop{\rm End}\nolimits}

\def\Gal{\mathop{\rm Gal}\nolimits}
\def\Ker{\mathop{\rm Ker}\nolimits}

\def\Im{\mathop{\rm Im}\nolimits}

\def\to{\rightarrow}

\def\normressym(#1,#2)_#3{\displaystyle\left({#1,#2\over#3}\right)}

\def\mod{\mathop{\rm mod.}\nolimits}
\def\pmod#1{\;(\mod#1)}

\newcount\refno 
\long\def\ref#1:#2<#3>{                                        
\global\advance\refno by1\par\noindent                              
\llap{[{\bf\number\refno}]\ }{#1} \pointir{\it #2} #3\goodbreak }

\def\citer#1(#2){[{\bf\number#1}\if#2\empty\relax\else,\ {#2}\fi]}

\newbox\bibbox
\setbox\bibbox\vbox{\bigbreak
\centerline{{\pc BIBLIOGRAPHIC} {\pc REFERENCES}}

\ref{\pc BELLA{\"I}CHE} (J.):
{\`A} propos d'un lemme de Ribet,
<Rend.\ Sem.\ Mat.\ Univ.\ Padova {\bf 109} (2003), 45--62.>
\newcount\bellaiche \global\bellaiche=\refno

\ref{\pc CASSELS} (J. W. S.):
Local fields,
<Cambridge University Press, Cambridge, 1986. xiv+360 pp.>
\newcount\cassels \global\cassels=\refno

\ref{\pc CASSELS} (J. W. S.) \& {\pc FR{\"O}HLICH} (A.):
Algebraic Number Theory,
<Academic Press, London, 1967, xviii+366 pp.>
\newcount\brighton \global\brighton=\refno

\ref{\pc COATES} (J.) \& {\pc SUJATHA} (R.):
The main conjecture,
<in Guwahati workshop on Iwasawa theory of totally real fields (Coates,
Dalawat, Saikia, Sujatha, Eds.), Ramanujan Mathematical Society, 2010.>
\newcount\coatessuja \global\coatessuja=\refno

\ref{\pc CURTIS} (C. W.) \& {\pc REINER} (I.):
Representation theory of finite groups and associative algebras.
<Interscience, New York, 1962, xiv+685 pp.>
\newcount\curtis \global\curtis=\refno

\ref{\pc DELIGNE} (P.):
Formes modulaires et repr{\'e}sentations $l$-adiques, 
<S{\'e}m.\ Bourbaki, Exp.\ 355, LNM {\bf 179},
Springer, Berlin, 1969, 139--172.> 
\newcount\deligne \global\deligne=\refno

\ref{\pc DELIGNE} (P.) \& {\pc RAPOPORT} (M.):
Les sch{\'e}mas de modules de courbes elliptiques. 
<in Modular functions of one variable II, LNM~349,
Springer, Berlin, 1973, 143--316.>
\newcount\dera \global\dera=\refno

\ref{\pc DELIGNE} (P.) \& {\pc SERRE} (J.-P.):
Formes modulaires de poids\/ $1$.
<Ann.\ sci.\ {\'E}cole norm.\ sup.\ (4) {\bf 7} (1975),  507--530.>
\newcount\delser \global\delser=\refno

\ref{\pc DIAMOND} (F.) \& {\pc SHURMAN} (J.):
A first course in modular forms.
<Springer-Verlag, New York, 2005, xvi+436 pp.>
\newcount\diamond \global\diamond=\refno

\ref{\pc FALTINGS} (G.) \& {\pc JORDAN} (B. W.):
Crystalline cohomology and ${\rm GL}(2,Q)$. 
<Israel J.\ Math.\ {\bf 90} (1995), no.~1-3, 1--66.>
\newcount\fajo \global\fajo=\refno

\ref{\pc FONTAINE} (J.-M.) \& {\pc LAFFAILLE} (G.):
Construction de repr{\'e}sentations $p$-adiques. 
<Ann.\ sci.\ {\'E}cole norm.\ sup.\ (4) {\bf 15} (1982), no.~4, 547--608.>
\newcount\fola \global\fola=\refno

\ref{\pc KHARE} (C.):
Notes on Ribet's converse to Herbrand. 
<Cyclotomic fields, Bhaskaracharya Pratishthana, Poona,
2000, 273--284.>
\newcount\khare \global\khare=\refno

\ref{\pc MAZUR} (B.):
How can we construc abelian Galois extensions of basic number fields.
<Preprint, 5 July 2010.>
\newcount\mazur \global\mazur=\refno

\ref{\pc RAYNAUD} (M.):
Sch{\'e}mas en groupes de type $(p,p,\ldots,p)$. 
<Bull.\ Soc.\ Math.\ France {\bf 102} (1974), 241--280.>
\newcount\raynaud \global\raynaud=\refno

\ref{\pc RIBET} (K. A.):
A modular construction of unramified $p$-extensions of\/ $\Q(\mu_p)$.
<Invent.\ Math.\ {\bf 34} (1976), no. 3, 151--162.>
\newcount\ribet \global\ribet=\refno

\ref{\pc RIBET} (K. A.):
Bernoulli numbers and ideal classes.
<Gaz.\ math.\ {\bf 118} (2008), 42--49.>
\newcount\gazette \global\gazette=\refno

\ref{\pc SAIKIA} (A.):
Ribet's construction of a suitable cusp eigenform.
<This volume.>
\newcount\anupam \global\anupam=\refno

\ref{\pc SERRE} (J.-P.):
Une interpr{\'e}tation des congruences relatives {\`a} la fonction $\tau$ de
Ramanujan. 
<S{\'e}m.\ DPP, Fasc.~1, Exp.~14, 17 pp. 
Secr{\'e}tariat math{\'e}matique, Paris, 1969.>
\newcount\inter \global\inter=\refno

\ref{\pc SERRE} (J.-P.):
Abelian $l$-adic representations and elliptic curves.
<W. A. Benjamin, New York-Amsterdam, 1968, xvi+177 pp.>
\newcount\serre \global\serre=\refno

\ref{\pc WASHINGTON} (L. C.):
Introduction to cyclotomic fields
<Springer-Verlag, New York, 1997. xiv+487 pp.>
\newcount\washington \global\washington=\refno

} 

\centerline{\bf Ribet's modular construction of unramified
  $p$-extensions of $\Q({}_p\mu)$}  
\medskip
\centerline{\it Notes for a series of lectures delivered at Gauhati}
\medskip
\centerline{Chandan Singh Dalawat}
\bigskip
\rightline{\it Ribet's ingenious method is the gateway to much later
    progress.}
\rightline{  --- Barry Mazur}

For every prime number $p$, we have the cyclotomic field $K=\Q(\zeta)$,
obtained by adjoining to $\Q$ a primitive $p$-th root $\zeta$ of~$1$~; it is a
galoisian extension of $\Q$, cyclic of degree~$p-1$.  The ideal class group $A$
of $K$ is a finite commutative group which measures the failure of the ring of
integers of $K$ to be factorial or principal.  We are interested in the
$\F_p$-space $C=A/A^p$. As $A$ and $C$ are trival for $p=2$, we assume that
the prime $p$ is {\it odd}.

In the middle of the XIX$^{\rm th}$ century, Kummer made a deep
connection between the group $C$ and the sequence of Bernoulli
numbers $(B_k)_{k\in\N}$ defined by 
$$
{T\over e^T-1}=\sum_{k\in\N}B_k{T^k\over k!},
$$
so that $B_0=1$, $B_1=-1/2$.  As the function $\displaystyle {T\over
  e^T-1}-1+{T\over 2}$ is even (invariant under $T\mapsto-T$), we
have $B_k=0$ for $k$ odd $>1$.  Here are the numerators $n_k$ and the
denominators $d_k$ of the first few Bernoulli numbers $B_k=n_k/d_k$ of even
index $k>0$~:
$$
\vbox{\halign{&\hfil$#$\quad\cr
\multispan{11}\hrulefill\cr
k=&2&4&6&8&10&12&14&16&18&20\cr 
\noalign{\vskip-5pt}
\multispan{11}\hrulefill\cr
n_k=&1&-1&1&-1&5&-691&7&-3617&43867&-174611\cr
d_k=&6&30&42&30&66&2730&6&510&798&330\cr
\multispan{11}\hrulefill.\cr}}
$$
Kummer's criterion says that {\it the\/ $\F_p$-space\/ $C$ is\/ $\neq\{1\}$
  if and only if\/ $p$ divides (the numerator of) one of the rational
  numbers\/ $B_2,B_4,\ldots,B_{p-3}$.}  Equivalently (see \S6),
$p\,|\,B_2B_4\ldots B_{p-3}$~; the sequence of such {\it irregular\/} primes
begins with $37$ (which divides $B_{32}=7709321041217/510$), $59|B_{44}$,
$67|B_{58}$, $101|B_{68}$, $103|B_{24}$, $\ldots$ and includes $283|B_{20}$,
$617|B_{20}$, $691|B_{12}$, $3617|B_{16}$, $43867|B_{18}$.

The $\F_p$-space $C$ carries an action of the group $\Delta=\Gal(K|\Q)$, for
which there is a canonical isomorphism $\chi:\Delta\to\F_p^\times$~; its
inverse sends $a\in\F_p^\times$ to the automorphism $\sigma_a$ of $K$ whose
restriction to the $p$-th roots of~$1$ is $\zeta\mapsto\zeta^a$.  Since every
character $\Delta\to\bar\F_p^\times$ is of the form $\chi^i$ for a unique
$i\in\Z/(p-1)\Z$, we have a direct sum decomposition
$C=\oplus_{i\in\Z/(p-1)\Z}C(\chi^i)$, where
$$
C(\chi^i)=\{c\in C\;|\;\sigma(c)=\chi^i(\sigma)c
\hbox{ for every }\sigma\in\Delta\}.
$$
In the early 30s, Herbrand refined one of the implications in
Kummer's criterion.  He showed that, {\it if\/ $k\in[2,p-3]$ is an even
integer such that\/ $C(\chi^{1-k})\neq0$, then\/ $p$ divides (the numerator
of)\/ $B_k$}.  This is an easy consequence of Stickelberger's theorem
\citer\washington(p.~101).

The main result in Ribet's paper is the converse of Herbrand's theorem.  In
other words, he showed that {\it if $p$ divides $B_k$ for some even integer
  $k\in[2,p-3]$, then $C(\chi^{1-k})\neq0$.}  

(We may allow $k=p-1$ in the statement ``\thinspace $p|B_k\Leftrightarrow
C(\chi^{1-k})\neq0$\thinspace'' of the Herbrand-Ribet theorem. We shall see in
\S6\ that in fact $p$ always divides the {\it denominator\/} of $B_{p-1}$.  We
must therefore show that $C(\chi)=0$.  This follows from Kummer theory and the
fact that the only degree-$p$ cyclic extension of $\Q$ unramified outside $p$
is the one contained in $\Q(\root p\of\zeta)$.) 

Let us end this introduction by mentioning how the Main Conjecture (for the
cyclotomic $\Z_p$-extension of $\Q$), as first proved by Mazur-Wiles, implies
a quantitative generalisation of the Herbrand-Ribet theorem~; see Ribet's own
account \citer\gazette(), also available on his homepage.

Instead of the group $C=A\otimes\Fp$, consider the $p$-part $D=A\otimes\Zp$ of
the ideal class group $A$ of $K=\Q(\zeta)$.  As for $C$, we have the
decomposition $D=\oplus_{i\in\Z/(p-1)\Z}D(\omega^i)$ as a $\Delta$-module,
where $\omega:\Delta\to\Z_p^\times$ is obtained from $\chi$ by composing with
the canonical (``Teichm{\"u}ller'') section $\F_p^\times\to\Z_p^\times$ of the
projection $\Z_p^\times\to\F_p^\times$~; every character
$\Delta\to\bar\Q_p^\times$ is of the form $\omega^i$ for a unique
$i\in\Z/(p-1)\Z$.  The Herbrand-Ribet theorem can be restated as the
equivalence ``\thinspace $p|B_k\Leftrightarrow
D(\omega^{1-k})\neq0$\thinspace'' for every even $k\in[2,p-3]$.

The Main Conjecture implies that, $k\in[2,p-3]$ being an even integer, {\it
  the group\/ $D(\omega^{1-k})$ has the same order as\/
  $\Z_p/L(\omega^{k-1},0)\Z_p$}, where, for any power $\varepsilon$ of
$\omega$,
$$
L(\varepsilon,0)={-1\over p}\sum_{a=1}^{p-1}\varepsilon(\sigma_a)a,\qquad
(\sigma_a:\zeta\mapsto\zeta^a).
$$
Because of the Kummer congruence $L(\omega^{k-1},0)\equiv-B_k/k\pmod p$
\citer\washington(p.~61), we have $p|B_k\Leftrightarrow p|L(\omega^{k-1},0)$,
so Herbrand-Ribet is indeed a consequence of Mazur-Wiles.  (There is as yet no
unconditional proof that the groups $D(\omega^{1-k})$ are cyclic.)

Although it is now possible, after the work of Thaine and Kolyvagin, to give a
more elementary proof \citer\washington(Chap.~15) of Herbrand-Ribet, and
indeed of Mazur-Wiles, the original approach using galoisian representations
attached to modular forms was a major source of inspiration for much
subsequent work and remains of paramount importance today.

Our aim here is to provide the reader with an overview of Ribet's modular
construction~(\S1--4)~; the numbering of the statements is the same as
in~\citer\ribet().  \S0~contains some results of a general nature which are
needed in the proof.  \S5~mentions an alternative approach which can now be
carried out.  Finally, \S6~gives Witt's proof of the von Staudt--Clausen
theorem on the denominators of~$B_k$.

\bigbreak
\centerline{\bf 0. Prerequisites}
\medskip

In this section we recall, without proof but with references to the
literature, a certain number of results which Ribet needs in order to carry
out his construction.  The reader may go directly to \S1--4 and consult
this~\S\ only if he needs to look up the statement of a result being used in
the proof.

{\it Equivariance of the reciprocity map}\pointir First, a special case of a
certain functorial property of the Artin reciprocity map.  Let\/ $M|K|\Q$ be
finite extensions such that $M|\Q$ is galoisian with group $G$, such that
$K|\Q$ is galoisian with group $\Delta$, and such that the group $H\subset G$
of automorphisms of $M|K$ is commutative.  Denote by ${\cal C}_K$ the
id{\`e}le class group of $K$, and by $\psi:{\cal C}_K\to H$ the reciprocity
map.  The group $\Delta$ acts on ${\cal C}_K$ as well as on $H$, by
conjugation on the latter.  We will need the fact that $\psi$ is equivariant
for these actions~:

\th THEOREM 0.1 \citer\brighton(p.~199)
\enonce
For every\/ $\sigma\in G$ (with image\/ $\bar\sigma\in\Delta$) and every\/
$x\in{\cal C}_K$, we have\/ $\psi(\bar\sigma x)=\sigma\psi(x)\sigma^{-1}$ in\/
$H$. 
\endth

{\it Locally algebraic characters\/}\pointir Let $\bar\Q$ be an algebraic
closure of $\Q$ and put $G_{\Q}=\Gal(\bar\Q|\Q)$.  We refer the reader to
\citer\serre(p.~III-2) for the definition of a {\it locally algebraic\/}
abelian representation $G_{\Q}\to GL_n(K)$ where $K$ is now a finite
extension of $\Q_p$.  Denote the $p$-adic cyclotomic character by
$\tilde\chi:G_{\Q}\to\Z_p^\times\to K^\times$.

\th PROPOSITION 0.2
\enonce
Every locally algebraic character\/ $\varphi:G_{\Q}\to K^\times$ is of the
form\/ $\varphi=\tilde\chi^n\omega$ for some $n\in\Z$ and some character
$\omega:G_{\Q}\to K^\times$ of finite order.
\endth

{\it The semisimplification of the reduction of a $p$-adic
  representation}\pointir The next result says that the reduction of a
$p$-adic representation may depend upon the choice of a stable lattice, but
the semisimplification does not. More precisely, let $K$ be a finite extension
of $\Q_p$, $\ogoth$ the ring of integers of $K$, $\pgoth$ the unique maximal
ideal of $\ogoth$, and $F=\ogoth/\pgoth$ the residue field.  Recall that an
$\ogoth$-lattice in a finite-dimensional $K$-space $V$ is a
sub-$\ogoth$-module generated by some $K$-basis of $V$.

Let $G$ be a profinite group --- a compact totally disconnected group
--- and $\rho:G\to GL(V)$ a (continuous) representation.  Among the
$\ogoth$-lattices in $V$, there is at least one which is $G$-stable,
for if $\Lambda\subset V$ is any $\ogoth$-lattice, the subgroup
$H=\rho^{-1}(GL(\Lambda))$ is open in $G$, hence the index ${\bf
  (}G:H{\bf )}$ is finite, and the $\ogoth$-lattice $\sum_{\sigma\in
  G/H}\sigma(\Lambda)$ is $G$-stable.

The choice of a $G$-stable $\ogoth$-lattice $\Lambda\subset V$ leads to a
(continuous) representation $\rho_\Lambda:G\to GL(\Lambda)\to
GL(\Lambda/\pgoth\Lambda)$ on a finite-dimensional $F$-space, called the
reduction of $\rho$ along $\Lambda$.  The theorem of Brauer-Nesbitt says that

\th PROPOSITION 0.3 \citer\curtis(p.~215) 
\enonce
The semisimplification\/ $\bar\rho$ of\/ $\rho_\Lambda$ depends only on\/
$\rho$, not on the choice of\/ $\Lambda$.  
\endth

{\it The galoisian representation associated to a weight-$2$ level-$p$ cuspidal
  eigenform}\pointir Next, we need some facts about modular forms~; for $N>0$
we have the congruence subgroups $\Gamma_1(N)\subset\Gamma_0(N)\subset{\bf
  SL}_2(\Z)$ whose elements are $\displaystyle\equiv\pmatrix{1&*\cr0&1\cr}$
and $\displaystyle\equiv\pmatrix{*&*\cr0&*\cr}\pmod N$ respectively. For
$k>1$, we denote by ${\goth S}_k(N)$ the space of cuspidal modular forms of
weight~$k$ for $\Gamma_1(N)$, and by ${\goth S}_k(N,\varepsilon)$ the subspace
on which $\Gamma_0(N)/\Gamma_1(N)=(\Z/N\Z)^\times$ acts via a given character
$\varepsilon:(\Z/N\Z)^\times\to\C^\times$, so that 
$$
{\goth S}_k(N)=\bigoplus_{\varepsilon:(\Z/N\Z)^\times\to\C^\times}
{\goth S}_k(N,\varepsilon).
$$
When $\varepsilon=1$ is trivial, ${\goth S}_k(N,1)$ is the space of weight-$k$ 
cuspidal modular forms for $\Gamma_0(N)$.

Let $f\in{\goth S}_2(p)$ be a normalised weight-$2$ cuspidal eigenform of
level~$p$, and write its $q$-expansion (with $q=e^{2i\pi\tau}$ and\/
$\tau\in{\goth H}$) as $f=\sum_{n>0}a_nq^n$.

\th PROPOSITION 0.4 \citer\diamond(p.~234)
\enonce
The subfield\/ $K_f=\Q(a_1,a_2,\ldots)$ of\/ $\C$ is of finite degree over\/
$\Q$.  Moreover, if\/ $\varepsilon:(\Z/p\Z)^\times\to\C^\times$ is the
character such that\/ $f\in{\goth S}_2(p,\varepsilon)$, then the field\/
$\Q(\varepsilon)$ generated by the values of\/ $\varepsilon$ is contained in
the number field\/ $K_f$. 
\endth

Recall the construction of the ($2$-dimensional) representations
$\rho_{f,\pgoth}:G_{\Q}\to GL(V_{f,\pgoth})$ attached to a normalised
eigenform $f\in{\goth S}_2(p,\varepsilon)$.  Here $\pgoth$ is any prime of
$K_f$ and $V_{f,\pgoth}$ is a certain vector space --- described below ---
over the completion $K_{f,\pgoth}$ of $K_f$ at the prime $\pgoth$.  Shimura
attaches to $f$ an abelian variety $A_f$ over $\Q$ (uniquely determined up to
$\Q$-isogeny) with the following properties.

The abelian variety $A_f$ has dimension $[K_f:\Q]$~; it is a quotient of the
jacobian $J_1(p)=\Jac X_1(p)$ of the modular curve $X_1(p)$~; and there is an
embedding of $\Q$-algebras $K_f\to\End_\Q(A_f)\otimes\Q$.
  
Let $T_p(A_f)=\lim_{n\in\N^*}\,{}_{p^n}A_f(\bar\Q)$ be the $p$-adic Tate
module of $A_f$ and $V_p(A_f)=T_p(A_f)\otimes_{\Z_p}\Q_p$, which is thus a
module over $K_f\otimes\Q_p$.

\th PROPOSITION 0.5 \citer\diamond(p.~390)
\enonce
The\/ $(K_f\otimes\Q_p)$-module\/ $V_p(A_f)$ is free of rank~$2$.
\endth
It follows that $V_{f,\pgoth}=V_p(A_f)\otimes_{K_f\otimes\Q_p} K_{f,\pgoth}$
is a $2$-dimensional vector space over $K_{f,\pgoth}$ with a continuous linear
action of $G_{\Q}$.  Denote this degree-$2$ representation by
$$
\rho_{f,\pgoth}:G_{\Q}\to GL(V_{f,\pgoth}).
$$

The abelian variety $A_f$ has good reduction at every prime $l\neq p$ of
$\Q$.  Further, the Eichler-Shimura relations give

\th PROPOSITION 0.6 \citer\diamond(p.~390)
\enonce
The representation\/ $\rho_{f,\pgoth}$ is unramified at every prime\/ $l\neq
p$ and the trace (resp.~determinant) of a Frobenius element $\Frob_l$ at $l$
is $a_l$ (resp.~$l\varepsilon(l)$) in $K_{f,\pgoth}$
(resp.~$K_{f,\pgoth}^\times$).  \endth

The final bit of information about $A_f$ that we need is the following theorem
of Deligne-Rapoport.

\th THEOREM 0.7 \citer\dera(p.~113)
\enonce
Suppose that\/ $\varepsilon$ is not trivial.  The abelian variety\/ $A_f$
acquires good reduction at the unique prime dividing $p$ in the fixed field\/
$\Q(\zeta)^+\subset\Q(\zeta)$ (where $\zeta^p=1$, $\zeta\neq1$) of the
involution\/ $\zeta\mapsto\zeta^{-1}$.
\endth

{\it Finite flat\/ group schemes}\pointir Let $E|\Qp$ be a finite extension,
with ring of integers $\ogoth_E$.  Ribet makes crucial use of the following
result of Raynaud.

\th THEOREM 0.8 \citer\raynaud(p.~268)
\enonce
Suppose that the ramification index of\/ $E|\Qp$ is $<p-1$.  Let\/ $G$ be
a finite flat commutative group scheme over $E$, killed by a power of\/ $p$.
There is at most one finite flat extension of\/ $G$ to\/ $\ogoth_E$.
\endth

\bigbreak
\centerline{\bf 1. The main theorem}
\medskip

Adopt the same notation as in the introduction.  In particular, $p$ is an odd
prime and $K=\Q(\zeta)$, with $\zeta^p=1$, $\zeta\neq1$.  We reduce the main
theorem (1.1) to the existence of a suitable finite extension $E$ (1.2) of
$\Q$, which in turn follows from the construction of a suitable representation
$r$ (1.3) of $\Gal(\bar\Q|\Q)$.  The existence of $r$ is established in \S4,
as a suitable reduction (\S2) of the representation associated to a certain
$f\in{\goth S}_2(p,\varepsilon)$ for a suitable $\varepsilon$ (\S3).

\th THEOREM 1.1
\enonce
Let\/ $k\in[2,p-3]$ be an even integer.  If\/ $p|B_k$, then
$C(\chi^{1-k})\neq0$.  
\endth
By the functorial property (0.1) of the Artin symbol, we are reduced to 
constructing a suitable everywhere-unramified extension $E|K$.

\th THEOREM 1.2
\enonce
Let\/ $k\in[2,p-3]$ be an even integer, and suppose that\/ $p|B_k$.
There exists a galoisian extension $E|\Q$ containing $K$ such that
\item {(a)} The extension $E|K$ is everywhere unramified.
\item {(b)} The group $H=\Gal(E|K)$ is a $p$-elementary commutative group
    ($\neq1$)~; equivalently, $H$ is an $\F_p$-space of dimension $>0$. 
\item {(c)} For every $\sigma$ in $G=\Gal(E|\Q)$ (with image
  $\bar\sigma$ in $\Delta$) and every $\tau$ in $H$,
$$
\sigma\tau\sigma^{-1}=\chi(\bar\sigma)^{1-k}.\tau
$$
\endth 
{\it Proof that\/} $(1.2)\Rightarrow(1.1)$.  Let $M$ be the maximal
$p$-elementary abelian extension of $K$ which is unramified at every place of
$K$~; by maximality, $M|\Q$ is galoisian.  The Artin map $\psi:{\cal
  C}_K\to\Gal(M|K)$, where ${\cal C}_K$ is the id{\`e}le class group of $K$,
is a homomorphism inducing an isomorphism $C\to\Gal(M|K)$.  Moreover, for
every $\sigma\in\Gal(M|\Q)$ (with image $\bar\sigma\in\Delta$) and every
$x\in{\cal C}_K$, we have $\psi(\bar\sigma x)=\sigma\psi(x)\sigma^{-1}$ (0.1).

In other words, the isomorphism (of $\F_p$-spaces) $C\to\Gal(M|K)$ induced by
$\psi$ is $\Delta$-equivariant, and hence induces, for every $i\in\Z/(p-1)\Z$,
an isomorphism $C(\chi^i)\to\Gal(M|K)(\chi^i)$.  If $(1.2)$ holds, the latter
group is $\neq0$ for $i=1-k$, because it has the quotient $H$.  It follows
that if $p|B_k$, then $C(\chi^{1-k})\neq0$, so $(1.1)$ holds (if (1.2) holds).

The next theorem shows how such an $E$ can be constructed from a certain kind
of representation of $G_{\Q}$, where $\bar\Q$ is an algebraic closure of
$\Q$ containing $K$.  Choose a decomposition group $D_p\subset G_{\Q}$ at
the prime~$p$, and let $\chi$ also stand for the composite
$G_{\Q}\to\Delta\to\F_p^\times$.

\th THEOREM 1.3
\enonce
Let\/ $k\in[2,p-3]$ be even, and suppose that\/ $p|B_k$.
There exists a finite extension\/ $F|\F_p$, and a (continuous)
representation\/ $ r:G_{\Q}\to GL_2(F)$, such that
\item {(i)} $ r$ is unramified at every prime\/ $l\neq p$.
\item {(ii)} $ r$ is an extension of\/ $\chi^{k-1}$ by the trivial
  character\/ $1$~:
$$
 r\sim\pmatrix{1&\gamma\cr 0&\chi^{k-1}\cr}
$$
for some map\/ $\gamma:G_{\Q}\to F$.
\item {(iii)} The order of\/ $\Im(r)$ is divisible by\/ $p$.
\item {(iv)} The order of\/ $\Im( r|_{D_p})$ is not divisible by\/
  $p$. 
\endth
{\it Proof that\/} $(1.3)\Rightarrow(1.2)$.  Let $ r$ be a representation as
in $(1.3)$; we have $\det(r)=\chi^{k-1}$, and there is a tower of fields
$E'|K'|\Q$, with $E'$ the field fixed by $\Ker(r)$ and $K'$ the field fixed by
$\Ker(\chi^{k-1})$~; we have $K'\subset K$, because
$\Ker(\chi)\subset\Ker(\chi^{k-1})$.  Put $H'=\Gal(E'|K')$, so that the
sequence $1\to H'\to\Im(r)\to\Im(\chi^{k-1})\to1$ is exact.  Notice that $H'$
is a $p$-elementary commutative group, for $\gamma$ gives rise to an injective
homomorphism $H'\to F$, and its order is $>1$, for the order of $\Im(r)$ is
divisible by $p$, by $(1.3)(iii)$.

Next, we remark that $E'|K'$ is everywhere unramified.  Indeed, for every
prime $l\neq p$, the extension $E'|\Q$ is unramified at $l$ by $(1.3)(i)$, and
hence so is $E'|K'$ at every prime $\goth l$ of $K'$ prime to $p$.  As for the
unique prime $\pgoth$ of $K'$ dividing $p$, the inertia subgroup at $\pgoth$
has order prime to $p$ by $(1.3)(iv)$, so $E'|K'$ is at worst tamely ramified,
and hence unramified, because $E'|K'$ is a $p$-extension.

Taking $E=E'K$, it is now clear that $(1.2)(a)$ and $(1.2)(b)$ hold for
$E|K|\Q$, because the analogous statements hold for $E'|K'|\Q$, as we have
just seen.  It remains to show that $(1.2)(c)$ holds as well.  This follows
from the fact that for every $\sigma,\tau\in G_{\Q}$, we have
$$
\pmatrix{1&\gamma(\sigma)\cr0&\chi^{k-1}(\sigma)\cr}
\pmatrix{1&\gamma(\tau)\cr0&1\cr}
\pmatrix{1&\gamma(\sigma)\cr0&\chi^{k-1}(\sigma)\cr}^{-1}
=\pmatrix{1&\chi^{1-k}(\sigma).\gamma(\tau)\cr0&1\cr},
$$
because of the identity
$\pmatrix{a&b\cr0&d}\pmatrix{1&x\cr0&1}\pmatrix{a&b\cr0&d}^{-1} =
\pmatrix{1&ad^{-1}.x\cr0&1\cr}$.  

Thus the implications $(1.3)\Rightarrow(1.2)\Rightarrow(1.1)$ are proved.  It
remains to find an $r$ as in (1.3) to complete the proof of the main theorem
(1.1).  

We shall take $r$ to be a suitable reduction (\S2) of the representation
$\rho_{f,\pgoth}$ (\S4) attached to a certain weight-$2$ level-$p$ cuspidal
eigenform $f$ (\S3).  A different approach, involving a weight-$k$ level-$1$
cuspidal eigenform $g$ instead of $f$, will be sketched in \S5.

\bigbreak
\centerline{\bf 2. Reductions of $p$-adic representations}
\medskip

Let $K|\Q_p$ be a finite extension, $\ogoth$ its ring of integers, $\pgoth$
the maximal ideal of $\ogoth$, and $F$ the residue field.  Let $V$ be a
$2$-dimensional $K$-space with a continuous linear action $\rho$ of a
profinite group $G$, and suppose that $\bar\rho$ (0.3) is reducible.  The two
characters $\varphi_1,\varphi_2:G\to F^\times$ such that
$\bar\rho\sim\varphi_1\oplus\varphi_2$ are uniquely determined by $\rho$ up to
order.  Also, for every $G$-stable $\ogoth$-lattice $\Lambda\subset V$, the
reduction $\rho_\Lambda$ is reducible, and its semisimplification is
isomorphic to $\varphi_1\oplus\varphi_2$ (0.3).  For a given $\Lambda$,
$\rho_\Lambda$ is isomorphic to one of
$$
\pmatrix{\varphi_1&*\cr0&\varphi_2\cr},\quad
\pmatrix{\varphi_1&0\cr*&\varphi_2\cr},
$$
and $\rho_\Lambda$ is semisimple if $\rho_\Lambda(G)$ has
order prime to $p$.  Ribet's lemma says that each of these two types
can be realised for some $\Lambda$, and, if $\rho$ is simple (but
$\bar\rho$ is not), then there is a $\rho_{\Lambda}$ which is not
semisimple and has any of the displayed types prescribed in advance.

\th PROPOSITION 2.1 (``Ribet's lemma'')
\enonce
Suppose that the degree-$2$ representation\/ $V$ of $G$ is simple but\/
$\bar\rho$ is not simple.  Then, for any ordering\/
$\varphi_1,\varphi_2$ of the two characters of which\/ $\bar\rho$ is the
direct sum, there is a\/ $G$-stable $\ogoth$-lattice\/ $\Lambda\subset V$
such that\/ $\rho_\Lambda\sim\pmatrix{\varphi_1&*\cr0&\varphi_2\cr}$,
as opposed to\/ $\pmatrix{\varphi_1&0\cr*&\varphi_2\cr}$, and such that\/ 
$\rho_\Lambda$ is not semisimple, i.e.,
$\rho_\Lambda\not\sim\varphi_1\oplus\varphi_2$. 
\endth
\def\\#1{\hbox{$\mit #1$}} Choosing a $G$-stable $\ogoth$-lattice
$\Lambda\subset V$ and fixing an $\ogoth$-basis of $\Lambda$ allows
one to view $\rho$ as a representation $G\to GL_2(\ogoth)$~; the
reduction of $\rho$ along $\Lambda$ is then the composite
$\rho_\Lambda:G\to GL_2(\ogoth)\to GL_2(F)$~; we call it the reduction
of $\rho(G)$.  For every $\\M\in GL_2(K)$ such that
$\\M\rho(G)\\M^{-1}\subset GL_2(\ogoth)$, we get another $G$-stable
$\ogoth$-lattice in $V$ together with an $\ogoth$-basis, giving
another reduction $G\to \\M\rho(G)\\M^{-1}\to GL_2(\ogoth)\to GL_2(F)$
of $\rho$. We will show that at least one of these reductions is of
the desired type $\pmatrix{\varphi_1&*\cr0&\varphi_2\cr}$, and that,
if all reductions of this type are semisimple, then $\rho$ cannot be
simple, which it is by hypothesis.

Let $\pi$ be an $\ogoth$-basis of $\pgoth$ (a uniformiser of $K$), and put
$\\P=\pmatrix{1&0\cr0&\pi\cr}$~; we have
$$\\P\pmatrix{a&\pi b\cr c&d\cr}\\P^{-1} 
 =\pmatrix{a&b\cr\pi c&d\cr}.
$$
This implies that there is at least one $G$-stable $\ogoth$-lattice in
$V$ along which the reduction of $\rho$ is of the desired type
$\pmatrix{\varphi_1&*\cr0&\varphi_2\cr}$, for if the reduction of
$\rho(G)$ is of the other type
$\pmatrix{\varphi_1&0\cr*&\varphi_2\cr}$, then the reduction of
$\\P\rho(G)\\P^{-1}$ is of the desired type, as the above formula
shows.

Suppose therefore that the reduction of $\rho(G)$ is of the desired
type $\pmatrix{\varphi_1&*\cr0&\varphi_2\cr}$~; it remains to show
that, among the reductions of $\rho$ of this type, there is at least
one which is not semisimple.  Suppose, if possible, that every
reduction of the type $\pmatrix{\varphi_1&*\cr0&\varphi_2\cr}$ is
semisimple.  We shall construct a sequence of matrices
$(\\M_i)_{i\in\N}$ in $GL_2(\ogoth)$, having a limit $\\M\in
GL_2(\ogoth)$, such that
$\\M\rho(G)\\M^{-1}\subset\pmatrix{\ogoth^\times&0\cr\pgoth&\ogoth^\times\cr}$,
showing that $\rho$ is not simple --- a contradiction.

Suppose that, for some
$i\in\N$, a matrix $\\M_i=\pmatrix{1&t_i\cr0&1\cr}$ in $GL_2(\ogoth)$
has been found such that
$$
\\M_i\rho(G)\\M_i^{-1}\subset
\pmatrix{\ogoth^\times&\pgoth^i\cr\pgoth&\ogoth^\times\cr}~;
$$
for $i=0$, take $\\M_0=\pmatrix{1&0\cr0&1\cr}$.  The above inclusion
can be rewritten $\\P^i\\M_i\rho(G)\\M_i^{-1}\\P^{-i}\subset
\pmatrix{\ogoth^\times&\ogoth\cr\pgoth^{i+1}&\ogoth^\times\cr}$, whose
reduction is of the desired type
$\pmatrix{\varphi_1&*\cr0&\varphi_2\cr}$, because the reduction of
$\rho(G)$ is of this type.  As all reductions of this type are
semisimple by hypothesis, there exists a matrix
$\\U_i=\pmatrix{1&u_i\cr0&1\cr}$ in $GL_2(\ogoth)$ which diagonalises
the reduction in the sense that
$$
\\U_i\\P^i\\M_i\rho(G)\\M_i^{-1}\\P^{-i}\\U_i^{-1}
 \subset\pmatrix{\ogoth^\times&\pgoth\cr\pgoth&\ogoth^\times\cr}.
$$
Denoting the entries of a matrix $()\in GL_2(\ogoth)$ by
$(a_{11}(),a_{12}();a_{21}(),a_{22}())$, one has
$a_{21}(\\U_i\\M\\U_i^{-1})=a_{21}(\\M)$.  This means that we actually
have
$$
\\U_i\\P^i\\M_i\rho(G)\\M_i^{-1}\\P^{-i}\\U_i^{-1}
 \subset\pmatrix{\ogoth^\times&\pgoth\cr\pgoth^{i+1}&\ogoth^\times\cr}, 
$$
for the $a_{21}$ of the matrix being conjugated by $\\U_i$ was in
$\pgoth^{i+1}$.  Conjugating by $\\P^{-i}$ and taking
$\\M_{i+1}=\\P^{-i}\\U_i\\P^i\\M_i=\pmatrix{1&t_i+u_i\pi^i\cr0&1\cr}$,
the above inclusion implies
$$
\\M_{i+1}\rho(G)\\M_{i+1}^{-1}\subset
 \pmatrix{\ogoth^\times&\pgoth^{i+1}\cr\pgoth&\ogoth^\times\cr},
$$
completing the construction of the sequence $(\\M_i)_{i\in\N}$ by
induction.  It is clear that this sequence converges in
$GL_2(\ogoth)$, and that the limit $\\M$ satisfies
$\\M\rho(G)\\M^{-1}\subset
 \pmatrix{\ogoth^\times&0\cr\pgoth&\ogoth^\times\cr}$,   
which contradicts the simplicity of $\rho$ and thereby proves that
$\rho$ has at least one reduction which has the desired type but which
is not semisimple.

For a generalisation of (2.1), see \citer\bellaiche().

\bigbreak
\centerline{\bf 3. A weight-$2$ level-$p$ cuspidal eigenform}
\medskip

Let $p$ be an odd prime, and let $\omega:\F_p^\times\to\Zp^\times$ be
the unique character (the Teichm{\"u}ller lift) such that
$\omega(d)\equiv d\pmod p$ for every $d\in\F_p^\times$.  The image of
$\omega$ is the subgroup ${}_{p-1}\mu$ or roots of~$1$ of order dividing
$p-1$. If a number field $K$ contains these roots of~$1$, and if we
fix a prime $\pgoth|p$ of $K$, then we may think of $\omega$ as taking
values in $K^\times$.

Let $k\in[2,p-3]$ be an even integer such that $p|B_k$, and put
$\varepsilon=\omega^{k-2}$.  See Anupam Saikia's contribution \citer\anupam()
to this volume for the proof of the following crucial result.

\th THEOREM 3.7
\enonce
Suppose that $p|B_k$.  There exist a normalised cuspidal eigenform\/
$f\in{\goth S}_2(p,\varepsilon)$, $f=\sum_{n>0}a_nq^n$, of weight\/~$2$,
level\/~$p$, and character\/~$\varepsilon$, and a prime\/ $\pgoth|p$ of the
number field\/ $K_f$ $(0.4)$ generated by the $q$-coefficients of $f$ such
that for every prime\/ $l\neq p$, the number\/ $a_l$ is\/ $\pgoth$-integral
and
$$
a_l\equiv1+l^{k-1}\equiv1+\varepsilon(l)l\pmod\pgoth.
$$
\endth 

Note that $K_f$ is a number field (0.4) and that $\varepsilon$ may be thought
of as taking values in $K_f^\times$, for $K_f$ contains $\Q(\varepsilon)$, the
number field generated by the values of $\varepsilon$ (0.4).

Now let $\rho_{f,\pgoth}$ be the representation attached to $f$ at $\pgoth$
(\S0)~; we show (\S4) that it has a reduction $r$ having the properties
$(i)$--$(iv)$ of (1.3).

\bigbreak
\centerline{\bf 4. The galoisian representation}
\medskip

Denoting the $p$-adic cyclotomic character by
$\tilde\chi:G_{\Q}\to\Z_p^\times$ (so that $\chi:\Delta\to\F_p^\times$ is
the reduction of $\tilde\chi$), and thinking of it as taking values in
$K_{f,\pgoth}^\times$, the equality
$\det(\rho_{f,\pgoth}(\Frob_l))=l\varepsilon(l)$ (0.6) can be rewritten as
$\det(\rho_{f,\pgoth})=\tilde\chi.\varepsilon$.

\th PROPOSITION 4.1
\enonce
The representation\/ $\rho_{f,\pgoth}$ is (absolutely) simple.
\endth
Suppose that it is not~; its semisimplification is then $\rho_1\oplus\rho_2$,
for some characters $\rho_i:G_{\Q}\to K_{f,\pgoth}^\times$.  It follows
(0.2) that each $\rho_i$ can be written as $\chi^{n_i}$ (for some $n_i\in\Z$)
on an open subgroup of the inertia subgroup at $p$, and therefore
$\rho_i=\varepsilon_i\chi^{n_i}$ for some characters of finite order
$\varepsilon_i:G_{\Q}\to K_{f,\pgoth}^\times$ unramified away from $p$.

Thinking of the $\varepsilon_i$ as characters of $(\Z/m_i\Z)^\times$
for some $m_i$, we have, for every prime $l\neq p$, the relations
$$
\varepsilon_1(l)l^{n_1}+\varepsilon_2(l)l^{n_2}=a_l,\quad
\varepsilon_1(l)\varepsilon_2(l)l^{n_1+n_2}=l\varepsilon(l).
$$
From the second equation we get $n_1+n_2=1$, and hence one of the
$n_i$ is $\ge1$ and the other $\le0$.  The first equation then gives
$|a_l|\ge l-1$ which, for $l\ge7$, contradicts the ``Riemann hypothesis''
$|a_l|\le2\sqrt l$, as proved by Eichler-Shimura.

Denote the ring of integer of $K_{f,\pgoth}$ by $\ogoth_{f,\pgoth}$.

\th PROPOSITION 4.2
\enonce
There exists a\/ $G_{\Q}$-stable\/
$\ogoth_{f,\pgoth}$-lattice\/ $\Lambda\subset V_{f,\pgoth}$ such that
$$
\rho_{f,\pgoth,\Lambda}\sim\pmatrix{1&*\cr0&\chi^{k-1}\cr},\quad
\rho_{f,\pgoth,\Lambda}\not\sim\pmatrix{1&0\cr0&\chi^{k-1}\cr}.
$$
\endth
We know that $\rho_{f,\pgoth}$ is irreducible (4.1).  The proposition would
follow from (2.1) once we show that
$\bar\rho_{f,\pgoth}\sim1\oplus\chi^{k-1}$, where $\bar\rho_{f,\pgoth}$ is the
semisimplification of $\rho_{f,\pgoth,\Lambda}$ (0.3).

Indeed, for every prime $l\neq p$, we know (0.6) that $\rho_{f,\pgoth}$, and
hence the reduction $\bar\rho_{f,\pgoth}$, is unramified at $l$.  Moreover
(0.6), the trace and the determinant of $\Frob_l\in\Im(\rho_{f,\pgoth})$ are
$a_l$ and $l\varepsilon(l)$ respectively.  By (3.7) and the definition of
$\varepsilon$, these numbers are congruent respectively to $1+l^{k-1}$ and
$l^{k-1}$ modulo $\pgoth$.  By the Chebotarev density theorem, the trace and
the determinant of $\bar\rho_{f,\pgoth}$ are $1+\chi^{k-1}$ and $\chi^{k-1}$
respectively, which are the same as the trace and the determinant of the
(semisimple) representation $1\oplus\chi^{k-1}$.  But two degree-$2$
semisimple representations having the same trace and the same determinant are
the same (0.3), hence $\bar\rho_{f,\pgoth}\sim1\oplus\chi^{k-1}$, which was to
be proved.

Now fix such a lattice $\Lambda$ and let $r=\rho_{f,{\pgoth},\Lambda}$ be the
reduction of $\rho_{f,\pgoth}$ along $\Lambda$.  We have seen that the
representation $r$ satisfies the conditions $(i)$--$(iii)$ of (1.3)~; it
remains to check that it satisfies $(iv)$ --- the restriction $r|_{D_p}$ of
$r$ to a decomposition group $D_p$ at $p$ is semisimple.

Denote by $F$ the residue field of $K_{f,\pgoth}$ and by $M$ the
(2-dimensional) $F$-space which carries the representation $r$.  Let $L$ be
the completion of $\Q(\zeta)^+$ at the unique place above $p$, and
$D=\Gal(\bar L|L)$ the decomposition group at the said place.  It is
sufficient to show that $r|_D$ (or the $D$-module $M$) is semisimple, for the
index $(D_p:D)$ is prime to~$p$, where we have identified $D_p$ with
$\Gal(\bar L|\Q_p)$.

\th PROPOSITION 4.3
\enonce
The\/ $D$-module\/ $M$ comes from a finite flat\/ $p$-elementary commutative
group scheme over the ring of integers\/ $\ogoth_L$ of\/ $L$.
\endth
There is an abelian $\Q$-variety $A$ which is $\Q$-isogenous to $A_f$
and such that $\ogoth_f\subset\End(A)$, because
$K_f\subset\End_\Q(A_f)\otimes\Q$.  The $D$-module $M$ is then the
``kernel'' ${}_\pgoth A$ of $\pgoth$, and hence contained in the
$D$-module ${}_pA$, the kernel of $p$.  As $A$ has good reduction over
$L$ (0.7), the $D$-module ${}_pA$ comes from a finite flat
$p$-elementary commutative group scheme over $\ogoth_L$, namely the
kernel ${}_p{\cal A}$ of $p$, where $\cal A$ --- the N{\'e}ron
model of $A$ --- is the abelian $\ogoth_L$-scheme whose generic fibre
is $A$.  It follows that the $D$-module $M={}_\pgoth A$ comes from the
the schematic closure $\cal M$ of $M$ in ${}_p\cal A$ (0.8).

\th THEOREM 4.4
\enonce
The\/ $D$-module\/ $M$ is semisimple.
\endth
This is the same as saying that the order of\/ $r(D)$ is prime
to\/~$p$, and will follow from the existence of two distinct $F$-lines ($X$
and $M^0$) in
$M$ which are $D$-stable.

We know that $r$, and hence $r|_D$, is isomorphic to
$\pmatrix{1&*\cr0&\chi^{k-1}\cr}$.  Let $X\subset M$ be an $F$-line on which
$D$ acts trivially, so that $D$ acts via $\chi^{k-1}$ on the $F$-line
$Y=M/X$~; in particular, $Y$ is ramified.

Let $\cal X$ be the schematic closure of $X$ in $\cal M$, so that the
$D$-module coming from $\cal X$ is $X$.  As the absolute ramification index
$(p-1)/2$ of $L$ is $<p-1$, the group scheme $\cal X$ is constant (0.8).  It
follows that $\cal M$ cannot be connected, for it has the {\'e}tale subgroup
scheme $\cal X$ (of order $>1$).

Now, the group scheme $\cal M$ is an $F$-space scheme.  Let $M^0\subset M$ be
the sub-$F[D]$-module coming from the largest connected subgroup scheme of
$\cal M$, so that the $F[D]$-module $M/M^0$ comes from the largest {\'e}tale
quotient of $\cal M$ and hence $M/M^0$ is unramified.

We have $M^0\neq M$, because $\cal M$ is not connected, as we have seen.  We
have $M^0\neq0$ because $M/M^0$ is unramified whereas $M$ is not (for it has
the quotient $Y$ which is ramified).  For the same reason, $M^0\neq X$,
because $M/M^0$ is unramified whereas $Y=M/X$ is ramified.  Thus $X$ and $M^0$
are two distinct $D$-stable $F$-lines in $M$, and hence the $D$-module $M$ is
semisimple, which was to be shown.

{\it Summary of the proof of the main theorem\/ $(1.1)$}\pointir The first
step (3.7) is to construct, if $p|B_k$, a normalised cuspidal eigenform
$f\in{\goth S}_2(p)$ such that for every prime $l\neq p$, the
$q$-coefficient $a_l(f)$ is $\equiv 1+l^{k-1}$ modulo a certain prime
$\pgoth|p$ of $K_f$.  A certain reduction $r$ (2.1) of the associated
$\pgoth$-adic representation $\rho_{f,\pgoth}$ gives an everywhere-unramified
$p$-elementary abelian extension $E|\Q(\zeta)$ which is galoisian over $\Q$ and
such that $\Delta=\Gal(\Q(\zeta)|\Q)$ operates on $H=\Gal(E|\Q(\zeta))$ via
the character $\chi^{1-k}:\Delta\to\F_p^\times$.  Hence $C(\chi^{1-k})\neq0$.

\bigbreak
\centerline{\bf 5. A weight-$k$ level-$1$ cuspidal eigenform}
\medskip

Ribet was aware at the time of writing \citer\ribet() of Serre's suggestion
\citer\inter() that instead of the weight-$2$ level-$p$ cuspidal eigenform
$f\in{\goth S}_2(p,\varepsilon)$ of (3.7), one could work with a suitable
weight-$k$ level-$1$ cuspidal eigenform $g\in{\goth S}_k(1)$.  The
construction of $g$ is more conceptual than that of $f$.  He could show that
for some prime ${\goth q}|p$ of the number field $K_g$, the associated
representation $\rho_{g,{\goth q}}$ has a reduction $r$ satisfying properties
$(i)$--$(iii)$ of (1.3).  

However, the techniques of the day were insufficient to prove that there is a
reduction satisfying $(1.3)(iv)$ as well.  But, as C.~Khare \citer\khare()
points out, it has now become possible to show that $\rho_{g,{\goth q}}$ has a
reduction $r$ satisfying all four requirements of (1.3), thanks to the work of
Faltings-Jordan \citer\fajo() and Fontaine-Laffaille \citer\fola().  We give
the broad outline.

\th PROPOSITION 5.1
\enonce
Suppose that $p|B_k$ for some even $k\in[2,p-3]$.  There exists a normalised
weight-$k$ level-$1$ cuspidal eigenform\/ $g=\sum_{n>0}a_n(g)q^n$ in\/ ${\goth
  S}_k(1)$ and a prime\/ ${\goth q}|p$ of the number field\/ $K_g$ generated
by the\/ $q$-coefficients of\/ $g$ such that for every prime\/ $l\neq p$, the
number\/ $a_l(g)$ is\/ ${\goth q}$-integral and
$$
a_l(g)\equiv1+l^{k-1}\pmod{\goth q}.
$$
\endth 
Recall that the Eisenstein series
$E_j=(-B_j/2j)+\sum_{n>0}\sigma_{j-1}(n)q^n$, where
$\sigma_i(n)=\sum_{d|n}d^i$ (so that $\sigma_i(l)=1+l^i$ for every prime~$l$),
is in the space ${\goth M}_j(1)$ of modular forms of weight $j$ (and
level~$1$) for every even integer $j>2$~; our $k$ is $>2$, for the numerator
of $B_2$ is $1$.  Notice that $a_l(E_k)=1+l^{k-1}$ for every prime $l$.  Also
recall that there is a unique normalised cuspidal eigenform $\Delta\in{\goth
  S}_{12}(1)$, with $q$-expansion
$$
q\prod_{n>0}(1-q^n)^{24}=\sum_{n>0}\tau(n)q^n.
$$
The graded $\C$-algebra of level-$1$ modular forms ${\goth M}(1)$ is the
polynomial ring $\C[E_4,E_6]$, and the graded ideal ${\goth S}(1)\subset{\goth
  M}(1)$ of level-$1$ cuspidal forms is principal, generated by $\Delta$
\citer\diamond(p.~88).

Following a suggestion of K.~Joshi, we therefore look for a $g$ of the type
$h=h_1\Delta+\cdots+h_n\Delta^n$ for a suitable $n>0$, where
$h_i=a_iE_4^{c_i}E_6^{d_i}$ ($a_i\in\Z$, $c_i,d_i\in\N$, $4c_i+6d_i+12i=k$).
Notice that every even integer $>2$ can be written as $4c+6d$ for some
$c,d\in\N$.

The existence of an $h\in{\goth S}_k(1)_{\Z}$ such that $h\equiv E_k\pmod p$
now follows from the fact that the constant terms $-B_4/8$ and $-B_6/12$ of
$E_4$, $E_6$ have numerators $1$, $-1$ respectively, and that if
$f_1,f_2\in{\goth M}_k(1)$ have the same first few (roughly $k/12$)
$q$-coefficients, then $f_1=f_2$, because the space ${\goth M}_k(1)$ is
finite-dimensional.  For example, $h=\Delta$ would do if $p=691$ and $k=12$. 

But we need a normalised {\it eigenform}, not just a cuspidal form. This is
furnished by the following lemma of Deligne-Serre, which implies the existence
of a number field $L'$, a prime ${\goth q}|p$ of $L'$, and a normalised
eigenform $g\in{\goth S}_k(1)_{L'}$ with ${\goth q}$-integral
$q$-coefficients, such that the congruence $a_l(g)\equiv 1+l^{k-1}\pmod{\goth
  q}$ holds for every prime $l\neq p$.

\th LEMMA 5.2 (Deligne-Serre \citer\delser(p.~522))
\enonce
Let\/ $A$ be a discrete valuation ring, $\goth m$ its maximal ideal, and\/ $L$
its field of fractions.  Let\/ $M$ be a free $A$-module of finite rank, and\/
$\goth h$ a set of commuting $A$-endomorphisms of\/ $M$.  Let\/ $f\in M/{\goth
  m}M$ be an eigenvector for every\/ $T\in{\goth h}$, of eigenvalue $a_T\in
A/{\goth m}$.  There exists a discrete valuation ring\/ $A'$ containing\/ $A$,
whose field of fractions is a finite extension of\/ $L$ and whose maximal
ideal\/ ${\goth m}'$ restricts to\/ $\goth m$ (${\goth m}'\cap A={\goth m}$),
and a common eigenvector\/ $f'\in M\otimes_AA'$ for\/ $\goth h$ such that the
eigenvalue\/ $a'_T\in A'$ of\/ $f'$ is\/ $\equiv a_T\pmod{{\goth m}'}$ for
every\/ $T\in{\goth h}$. 
\endth

See Anupam Saikia's contribution to this volume \citer\anupam(4.2) for a
proof.

To find a $g$ as in (5.1), apply this lemma with $A=\Z_{(p)}$, $M={\goth
  S}_k(1)_{\Z_{(p)}}$, ${\goth h}$ the set of Hecke operators $T_l$ for
primes~$l$, and $f=\bar h$, the reduction of $h$ modulo~$p$, which is an
eigenform in $M/pM={\goth S}_k(1)_{\F_p}$, the eigenvalue for $T_l$ being
$1+\bar l^{k-1}\in\F_p$ for every $l\neq p$, because $g\equiv E_k$.  Taking
$g$ to be the $f'$ furnished by (5.2) completes the proof of (5.1).

\smallskip

Recall that to $g$ and ${\goth q}$ as in (5.1) Deligne \citer\deligne()
attaches a degree-$2$ representation $\rho_{g,{\goth q}}$ of $G_{\Q}$ over
the completion $K_{g,{\goth q}}$ of $K_g=\Q(a_n(g))_{n>0}$ at ${\goth q}$~; it
is unramified at every prime $l\neq p$~; for these $l$, the trace of
$\rho_{g,{\goth q}}(\Frob_l)$ is $a_l(g)$, the determinant $l^{k-1}$.

\th PROPOSITION 5.3
\enonce
The representation\/ $\rho_{g,{\goth q}}$ is (absolutely) simple.
\endth

The proof is the same as for the representation $\rho_{f,\pgoth}$ in (4.1),
using the Ramanujan-Peterson estimate ``$|a_l(g)|\le2\sqrt{l^{k-1}}$ for every
prime $l\neq p$'', as proved by Deligne \citer\deligne(), instead of the
earlier ``$|a_l(f)|\le2\sqrt l$''.

\th PROPOSITION 5.4
\enonce
The representation\/ $\rho_{g,{\goth q}}$ has a reduction $r$ such that
$$
r\sim\pmatrix{1&*\cr0&\chi^{k-1}\cr},\quad
r\not\sim\pmatrix{1&0\cr0&\chi^{k-1}\cr}.
$$
\endth
This follows from (2.1) much in the same way as (4.2) does, because
$\bar\rho_{g,{\goth q}}\sim1\oplus\chi^{k-1}$, since the trace of every
reduction is $\equiv1+\chi^{k-1}\pmod{\goth q}$ and the determinant is
$\equiv\chi^{k-1}$.

\th PROPOSITION 5.5
\enonce
We may further require that the restriction\/ $r|_{G_p}$ of\/ $r$ to the
decomposition group\/ $G_p=\Gal(\bar\Qp|\Qp)$ at\/ $p$ is\/ semisimple~; in
fact, the order of\/ $r(G_p)$ is prime to\/ $p$. 
\endth
Let $F$ be the residue field of $K_{g,{\goth q}}$, so that $r|_{G_p}$ is a
degree-$2$ $F$-representation of $G_p$.  As for the $D$-module
$M$ (4.4), it suffices to show that the $F$-plane which carries the
representation $r|_{G_p}$ is the direct sum of two (distinct) $G_p$-stable
$F$-lines, to conclude that $r|_{G_p}$ is semisimple.

We have seen that there is a $G_p$-stable $F$-line on which $G_p$ acts
trivally~; it acts via $\chi^{k-1}$ on the quotient (5.4).  I am grateful to
Christophe Breuil for providing the following proof that there is also an
$G_p$-stable $F$-line on which $G_p$ acts via $\chi^{k-1}$, the action on the
quotient being trivial.

The basic idea is that $r|_{G_p}$ is (a reduction of) a {\it crystalline\/}
representation, by Faltings-Jordan \citer\fajo() and, as the weight $k$ is
$<p-1$, Fontaine-Laffaille \citer\fola() implies that the sequence
$0\to\chi^0\to r|_{G_p}\to\chi^{k-1}\to0$ must split.  We omit the details.

In other words, $(1.3)(i)$--$(iv)$ hold for $r$. This completes our sketch of
the proof of (1.3) using the level-$1$ weight-$k$ cuspidal eigenform $g$
instead of the level-$p$ weight-$2$ cuspidal eigenform $f$ of \S3.  As we have
seen in \S1, (1.3)$\Rightarrow$(1.1).

\medbreak

There is yet another approach to finding, when $p|B_k$, a cuspidal eigenform
$g\in{\goth S}_k(1)$, a number field $K$ containing all $a_n(g)$, and a prime
${\goth q}|p$ where all $a_n(g)$ are locally integral, such that $g\equiv
E_k\pmod{\goth q}$.  I'm grateful to Barry Mazur for pointing out this method
(briefly sketched by Khare), which offers certain advantages.

Denote by ${\goth H}\subset\C$ the open subset consisting of those $\tau=x+iy$
($i^2=-1$) for which $y>0$, and by ${\goth H}^*\subset\P_1(\C)$ the union
${\goth H}\cup\P_1(\Q)$.  The group $SL_2(\Z)$ acts on the spaces $\goth H$
and ${\goth H}^*$~; the $j$-invariant identifies the quotient ${\goth
  H}/SL_2(\Z)$ (resp.~${\goth H}^*\!/SL_2(\Z)$) with the affine
(resp.~projective) line and furnishes a canonical model $Y$ (resp.~$X$) of
this algebraic curve over $\Q$ (and indeed over $\Z$), so that for example
$X(\C)={\goth H}^*\!/SL_2(\Z)$ as analytic $\C$-curves.

For any commutative ring $R$, the space ${\goth S}_k(1)_R$ of cuspidal modular
forms of weight $k=2k'$, level~$1$, and with $q$-coefficients in $R$, can be
interpreted as the space $H^0(X_R,\Omega^{k'})$ of global sections of the line
bundle $\Omega^{k'}$ on the $R$-curve $X_R$, where $\Omega^1$ is the canonical
bundle, by sending the cuspidal form $f$ to the differential form
$f.(dq/q)^{k'}$.

Now, although the modular form $E_k\in{\goth M}_k(1)$ is not cuspidal, our
hypothesis $p|B_k$ implies that its constant term $-B_k/2k$ is ($p$-integral
and) $\equiv0\pmod p$.  Therefore the mod-$p$ reduction $\bar E_k$ of $E_k$ is
in ${\goth S}_k(1)_{\Fp}$, giving rise to the global section $\bar
E_k.(dq/q)^{k'}\in H^0(X_{\Fp},\Omega^{k'})$.  At the same time, the reduction
map $H^0(X_{\Zp}, \Omega^{k'})\to H^0(X_{\Fp},\Omega^{k'})$ is surjective, at
least if $p\neq2,3$ (which is the case here~: the odd prime $p$ is irregular
by hypothesis), so there is a cuspidal modular form $h\in{\goth S}_k(1)_{\Zp}$
such that $h\equiv\bar E_k\pmod p$.

But $E_k$ is an eigenform, the eigenvalue for $T_l$ being $1+l^{k-1}$ for all
primes~$l$, including $l=p$.  We can thus invoke the Deligne-Serre lemma (5.2)
to get a number field $K$, a cuspidal {\it eigenform\/} $g\in{\goth S}_k(1)$
whose $q$-coefficients are in $K$, and a place ${\goth q}|p$ of $K$ such that
the $a_n(g)$ are $\goth q$-integral and $a_l(g)\equiv 1+l^{k-1}\pmod{\goth q}$
for every prime $l$, thereby establishing (5.1).  Note in particular that
$a_p(g)\equiv1\pmod{\goth q}$.

Let $\rho_{g,{\goth q}}$ be the degree-$2$ representation
$G_{\Q}\to\Aut_{K_{\goth q}} (V)$ associated to $g$ and $\goth q$.  Ribet's
lemma (2.1) then provides a $G_{\Q}$-stable lattice $\Lambda\subset V$ such
that the reduction $r=\rho_{g,{\goth q},\Lambda}:G_{\Q}\to\Aut_{F}
(\Lambda\otimes_{{\goth o}_{\goth q}}F)$ of $\rho_{g,{\goth q}}$ along
$\Lambda$ sits in an unsplittable exact sequence
$$
0\to\chi^0\to r\to\chi^{k-1}\to0,
$$
just as in (4.2) and (5.4).  In other words, the $G_{\Q}$-module
$\Lambda\otimes_{{\goth o}_{\goth q}}F$ is reducible but indecomposable and
contains a $G_{\Q}$-stable $F$-line with the trivial action, the action on the
quotient being via $\chi^{k-1}$, so that its semisimplification is
$\sim1\oplus\chi^{k-1}$.  Here, ${\goth o}_{\goth q}$ is the ring of integers
of the completion $K_{\goth q}$ of $K$ at $\goth q$, and $F$ (a finite
extension of $\Fp$) is the residue field of ${\goth o}_{\goth q}$.

It remains to show that $\Lambda$ can be so chosen that moreover the
restriction $r|_{G_p}$ of $r=\rho_{g,{\goth q},\Lambda}$ to the decomposition
group $G_p=\Gal(\bar\Qp|\Qp)$ is semisimple.  This follows from the fact that
$g$ is $p$-ordinary at $\goth q$.  Indeed, as the character $\chi^{k-1}$ is
ramified, $\Lambda$ can be so chosen that there is also an exact sequence
$$
0\to\chi^{k-1}\to r|_{G_p}\to\chi^0\to0,
$$
which, in conjunction with the previous exact sequence (restricted to
$G_p$) implies that $r|_{G_p}$ is semisimple.  Mazur calls it the Ribet
``wrench'' in \citer\mazur(),~\S16~; we omit the details.  See also
\citer\coatessuja(),~\S6.

\bigbreak
\centerline{\bf 6. The denominators of Bernoulli numbers}
\medskip

Let's conclude with Witt's unpublished proof, as presented in
\citer\cassels(), about the prime divisors of the denominators of $B_k$.  It
implies in particular that the numerator of $B_{p-1}$ is never divisible by
$p$.  

\th THEOREM 6.1 (von Staudt--Clausen, 1840)  
\enonce 
Let $k>0$ be an even integer, and let $l$ run through the primes.  Then the
number 
$$
W_k=B_k+\sum_{l-1\,|\,k}{1\over l}\leqno{(1)}
$$
is always an integer.  For example, $\displaystyle
W_{12}=B_{12}+{1\over 2}+{1\over3}+{1\over5}+{1\over7}+{1\over 13}=1$.  
\endth

{\eightpoint The British analyst Hardy says in his {\it Twelve lectures\/}
  (p.~11) that this theorem was rediscovered by Ramanujan ``\thinspace at a
  time of his life when he had hardly formed any definite concept of
  proof\thinspace''.}

{\it Proof\/} (Witt)~: The idea is to show that $W_k$ is a $p$-adic integer
for every prime $p$.  More precisely, we show that $B_k+p^{-1}$ (resp.~$B_k$)
is a $p$-adic integer if $p-1\,|\,k$ (resp.~if not).

For an integer $n>0$, let $S_k(n)=0^k+1^k+2^k+\cdots+(n-1)^k$.  Comparing the
coefficients on the two sides of
$$
1+e^T+e^{2T}+\cdots+e^{(n-1)T}={e^{nT}-1\over T}{T\over e^T-1},
$$
we get $\displaystyle S_k(n) =\sum_{m\in[0,k]}{k\choose m}{B_m\over
  k+1-m}n^{k+1-m}$.  To recover $B_k$ from the $S_k(n)$, it is tempting to
take the limit $\displaystyle\lim_{n\to0}S_k(n)/n$, which doesn't make sense
in the archimedean world.  If, however, we make $n$ run through the powers
$p^s$ of a fixed prime $p$, then, $p$-adically, $p^s\to0$ as $s\to+\infty$,
and
$$
\lim_{s\to+\infty}S_k(p^s)/p^s=B_k.\leqno{(2)}
$$
Let us compare $S_k(p^{s+1})/p^{s+1}$ with $S_k(p^s)/p^s$.  Every
$j\in[0,p^{s+1}[$ can be uniquely written as $j=up^s+v$, where $u\in[0,p[$ and
$v\in[0,p^s[$.  Now,
$$
\eqalign{
S_k(p^{s+1})
&=\sum_{j\in[0,p^{s+1}[}j^k=\sum_{u\in[0,p[}\sum_{v\in[0,p^s[}(up^s+v)^k\cr
&\equiv p\left(\sum_v v^k\right)
 +kp^s\left(\sum_u u\sum_v v^{k-1}\right)\pmod{p^{2s}}\cr 
}$$
by the binomial theorem.  As $\sum_{v}v^k=S_k(p^s)$ and
$2\sum_uu=p(p-1)\equiv0\pmod p$, we get
$$
S_k(p^{s+1})\equiv pS_k(p^s)\pmod{p^{s+1}},
$$
where, for $p=2$, the fact that $k$ is even has been used.  Dividing
throughout by $p^{s+1}$, this can be expressed by saying that
$$
{S_k(p^{s+1})\over p^{s+1}}-{S_k(p^s)\over p^s}\in\Z_{(p)}
$$
is a $p$-adic integer, and therefore
$$
{S_k(p^r)\over p^r}-{S_k(p^s)\over p^s}\in\Z_{(p)}
$$
for any two integers $r>0$, $s>0$, since $\Z_{(p)}$ is a subring of $\Q$.
Fixing $s=1$ and letting $r\to+\infty$, we see that $B_k-S_k(p)/p\in\Z_{(p)}$,
in view of $(2)$.  We need

\th LEMMA 6.2
\enonce
$\displaystyle
S_k(p)=\sum_{j\in[1,p[}j^k\equiv\cases
{-1\pmod p&if $p-1\,|\,k$\cr
\phantom{-}0\pmod p&otherwise.\cr}
$
\endth
This is clear if $p-1\,|\,k$.  Suppose not, and let $g$ be a generator of
$(\Z/p\Z)^\times$.   We have $g^k-1\not\equiv0$, whereas
$$
(g^k-1)\left(\sum_{j\in[1,p[}j^k\right)
\equiv (g^k-1)\left(\sum_{t\in[0,p-1[}g^{tk}\right)
\equiv g^{(p-1)k}-1\equiv0.
$$
It follows that $B_k+p^{-1}\in\Z_{(p)}$ if $p-1\,|\,k$ and $B_k\in\Z_{(p)}$
otherwise.  In either case, the number $W_k$ (1), which can be written as
$$
W_k=\cases{
(B_k+p^{-1})+\sum_{l\neq p}l^{-1}&if $p-1\,|\,k$\cr
(B_k)+\sum_l l^{-1}&otherwise,\cr
}$$ 
(where $l$ runs through the primes for which $l-1\,|\,k$) turns out to be a
$p$-adic integer for every prime $p$.  Hence $W_k\in\Z$, as claimed.  

John Coates remarked at the workshop that the analogue of (6.1) for a totally
real number field $F$ (other than $\Q$) is an open problem~; even a weak
analogue would imply Leopoldt's conjecture for $F$.

\bigbreak
\unvbox\bibbox 
\bigbreak\vskip2cm
{\obeylines\parskip=0pt\parindent=0pt
Chandan Singh Dalawat
Harish-Chandra Research Institute
Chhatnag Road, Jhunsi
{\pc ALLAHABAD} 211\thinspace019, India
\vskip5pt
\tt dalawat@gmail.com
\vskip5pt
}

\bye